\documentclass[final]{siamonline1116} 

\usepackage{amsfonts,amssymb,amsmath,overpic}
\usepackage{algorithm}
\usepackage{algpseudocode}
\algnewcommand\algorithmicinput{\textbf{Input:}}
\algnewcommand\Input{\item[\algorithmicinput]}
\usepackage{array}
\newcolumntype{t}{>{\ttfamily}c}

\newcommand*\R{\mathbb{R}}

\newcommand{\ale}[2][red]{\emph{\textcolor{#1}{#2}}}
\graphicspath{{figures/}}

\newsiamremark{remark}{Remark}

\crefname{algorithm}{algorithm}{algorithms}
\crefname{remark}{remark}{remarks}

\title{Data-driven identification of parametric partial differential equations
  \thanks{Submitted 6/2/2018
    \funding{S. Rudy, S. L. Brunton and J. N. Kutz acknowledge support from the Defense Advanced Research Projects Agency (DARPA HR0011-16-C-0016).    J. N. Kutz acknowledges support from the Air Force Office of Scientific Research
      (AFOSR) grant FA9550-15-1-0385.  S. L. Brunton acknowledges support from the Air Force Office of Scientific Research (AFOSR FA9550-18-1-0200). }}
      }

\author{Samuel Rudy\thanks{Department of Applied Mathematics,
    University of Washington, Seattle, WA
    (\email{shrudy@uw.edu}, \email{kutz@uw.edu}).}
   \and
   Alessandro Alla\thanks{Department of Mathematics, PUC-Rio, Rua Marques de Sao Vicente, 225, Rio de Janeiro,  22451-900, Brasil (\email{alla@mat.puc-rio.br})}
  \and
  Steven L. Brunton\thanks{Department of Mechanical Engineering,
    University of Washington, Seattle, WA (\email{sbrunton@uw.edu})} 
  \and
  J. Nathan Kutz\footnotemark[2]
}

\begin{document}

\maketitle
\begin{abstract}
In this work we present a data-driven method for the discovery of parametric partial differential equations (PDEs), thus allowing one to disambiguate between the underlying evolution equations and their parametric dependencies.  Group sparsity is used to ensure parsimonious representations of observed dynamics in the form of a parametric PDE, while also allowing the coefficients to have arbitrary time series, or spatial dependence.  This work builds on previous methods for the identification of constant coefficient PDEs, expanding the field to include a new class of equations which until now have eluded machine learning based identification methods.  We show that group sequentially thresholded ridge regression outperforms group LASSO in identifying the fewest terms in the PDE along with their parametric dependency.  The method is demonstrated on four canonical models with and without the introduction of noise.
\end{abstract}




\section{Introduction}

The extraction of physical laws from experimental data, often in the form of differential and partial differential equations, may be critical to science and engineering applications where governing equations are unknown.  Time-series data collected from experiments, or as the macroscopic aggregate of small scale behavior, often obeys unknown governing equations that are parametrized by time-evolving parameters $\mu(t)$.  In some cases, it may be possible to derive physical laws from first principals using data as well as some knowledge of the system, but there are many cases where this is elusive such as the large scale networked dynamical systems of the power grid and the brain, or chemical kinetics of, for example, the  Belousov-Zhabotinsky reaction.  Recently, there has been a substantial research effort towards automating the process of data-driven model discovery in order to identify interpretable expressions for the dynamics in the form of ordinary and partial differential equations (ODEs and PDEs):
\begin{equation}
u_t= N(u, u_x, u_{xx}, \dots ,\mu(t)),\quad t\in [0,T],
\label{eq:overview}
\end{equation}
where $N(\cdot)$ characterizes the evolution of the system and its parametric dependencies through the parameter $\mu(t):[0,T]\rightarrow \R$.  Although a number of automated discovery techniques have been developed for discovering the right hand side of (\ref{eq:overview}) with constant parameters $\mu(t)=\mu_0$, none have demonstrated the capability to infer the governing equations when the parametrization $\mu(t)$ has explicit time dependence.  By imposing group sparsity techniques, we develop a mathematical architecture that allows us to explicitly disambiguate the dynamical evolution (\ref{eq:overview}) from its parametric dependencies $\mu(t)$.  This is a critical innovation in model discovery as most realistic systems do indeed have time-dependent parametric dependencies that must be concurrently extracted during the model discovery process.

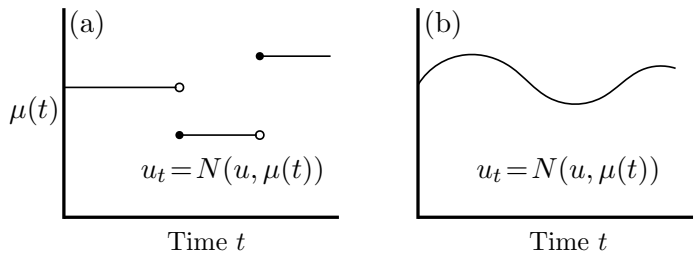
\begin{figure}[t]
\centering
\begin{overpic}[width=0.56\textwidth]{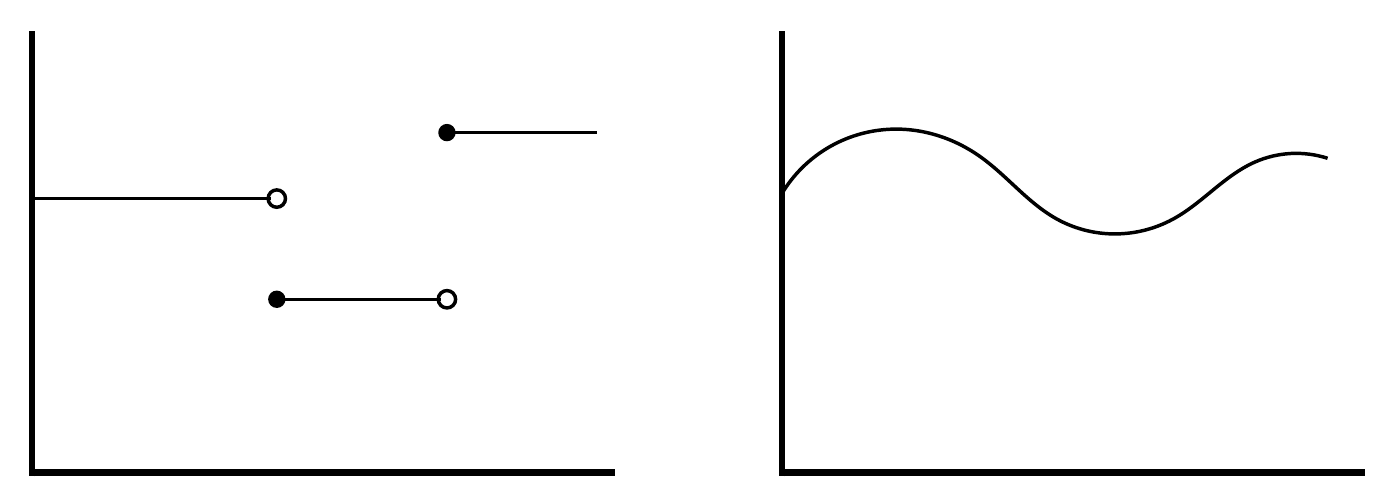}
\put(-6,17){${\mu}(t)$}
\put(3,30){(a)}
\put(57,30){(b)}
\put(18,-3){\small Time $t$}
\put(73,-3){\small Time $t$}
\put(65,8){$u_t \!=\! N({u},{\mu}(t))$}
\put(14,8){$u_t \!=\! N({u},{\mu}(t))$}
\end{overpic}
%
\caption{Two prototypical scenarios for parametric model discovery:  (a) Parameters $\mu(t)$ are piece wise constant and change at fixed points in time, (b) The underlying PDE model $u_t \!=\! N({u},{\mu}(t))$ depends on continuously varying parameter  ${\mu}(t)$.}
\label{fig:param}
\vspace{-.25in}
\end{figure}

Figure~\ref{fig:param} demonstrates two prototypical parametric dependencies:  (a) a PDE model whose constant parameters change at fixed points in time and (b) a PDE model that depends continuously on the parameter  $\mu(t)$.  Our proposed model discovery method provides a principled approach to efficiently handle these two parametric cases, thus advancing the field of model discovery by disambiguating the dynamics from parametric time dependence.
Even if the governing equations are known, parametric dependencies in PDEs complicate numerical simulation schemes and challenge one's ability to produce accurate future state predictions.  
Characterizing parametric dependence is also an critical task for model reduction in both time-independent~\cite{quarteroni2015reduced,hesthaven2016certified} and time-dependent PDEs~\cite{benner2015survey,benner2017model}.
Thus the ability to explicitly extract the parametric dependence of a spatio-temporal system is necessary for accurate, quantitative characterization of PDEs.

More broadly, system identification using machine learning methods has emerged as a viable alternative to expert knowledge and first principles derivations.  It is important to separate the field of system identification into two distinct categories:  (i) methods that accurately reflect observed dynamics using black box functions (e.g. neural networks), and (ii) methods that recover closed form and interpretable expressions for the dynamics in the form of ordinary and partial differential equations (ODEs and PDEs).   This duality reflects the two cultures narrative of machine learning and classical statistics made popular by Leo Breiman~\cite{breiman2001statistical}.  On one hand, the research may assume a specific model for the data with known mechanism, while on the other the research is interested in algorithmic models that, while not necessarily reflecting the true mechanism, are accurate in prediction.  While several recent works have made progress from both viewpoints, we focus on the former.  The terms in a differential equation often have physical interpretations and motivation, e.g. diffusion and advection are ubiquitous in many physical systems and are characterized by prototypical expressions in a PDE.  For systems where first principal derivations prove intractable, we may gain insights into the underlying physics of the system based on the terms in the identified PDE.  Further, we view the process of extracting closed form equations as being more generalizable than fitting black box models to a specific dataset.  Specifically, altering initial conditions will not change governing equations, but may break a machine learned black box solver.

Research towards the automated inference of dynamical systems from data is not new \cite{Crutchfield1987cs}.  Methods for extracting linear systems from time series data include 
the eigensystem realization algorithm (ERA) \cite{juang1985eigensystem} and Dynamic Mode Decomposition (DMD) \cite{Rowley2009jfm,schmid2010dynamic,tu2013dynamic,Brunton2015jcd,kutz2016dynamic,askham2018variable}.  Identification of nonlinear systems has, until very recently, relied on black box methods.  These include NARMAX \cite{chen1989representations}, neural networks \cite{gonzalez1998identification}, equation free methods \cite{Kevrekidis2003cms,kevrekidis2004equation,kevrekidis2009equation}, and Laplacian spectral analysis \cite{Giannakis2012pnas}.  There has also been considerable recent work towards data-driven approximation of the Koopman operator~\cite{Mezic2005nd,Budivsic2012chaos,Mezic2013arfm} via extensions of DMD \cite{williams2015data}, diffusion maps \cite{giannakis2017data}, delay-coordinates \cite{brunton2017chaos, Arbabi2016arxiv, das2017delay} and neural networks \cite{Yeung2017arxiv,Takeishi2017nips,Wehmeyer2017arxiv,Mardt2017arxiv,Otto2017arxiv,lusch2017deep}.

The use of genetic algorithms for nonlinear system identification \cite{Bongard2007pnas, Schmidt2009science} allowed for the derivation of physical laws in the form of ordinary differential equations.  Genetic algorithms are highly effective in learning complex functional forms but are slower computationally than simple regression.  Sparsity promoting methods have been used previously in dynamical systems~\cite{Schaeffer2013pnas, mackey2014compressive, caflisch2013pdes}, and sparse regression has been leveraged to identify parsimonious ordinary differential equations from a large library of allowed functional forms~\cite{Brunton2016,chen2017network}.  Much work building on the sparse regression framework has followed and includes inferring rational functions \cite{Mangan2016}, the use of information criteria for model validation~\cite{mangan2017model}, constrained regression for conservation laws \cite{loiseau2018constrained}, model discovery with highly corrupt data \cite{tran2017exact}, the learning of bifurcation parameters \cite{schaeffer2017learning_group}, stochastic dynamics \cite{boninsegna2017sparse}, weak forms of the observed dynamics \cite{schaeffer2017sparse}, and regression with small amounts of data \cite{schaeffer2017extracting, kaiser2017sparse}.  In contrast to sparse regression, a neural network based approach was proposed to identify ordinary differential equations, sacrificing some interpretability for a richer class of allowed functional forms \cite{raissi2018multistep}.

Sparse regression based methods for PDEs were first used in \cite{rudy2017data,schaeffer2017learning}.  These methods were demonstrated on a large class of PDEs and have the benefit of being highly interpretable, but struggle with numerical differentiation of noisy data.  In Rudy {\em et al.}~\cite{rudy2017data} the noise was addressed by testing with only a small degree of noise (large SNR), while in Schaeffer {\em et al.}~\cite{schaeffer2017learning} noise was added to the time derivative after it was computed from clean data.  Alternatively, Gaussian processes were used to determine linear PDEs \cite{raissi2017machine} and nonlinear PDEs known up to a set of coefficients \cite{raissi2017hidden}.  Using Gaussian process regression requires less data than sparse regression and naturally manages noise, but the method is only applicable to PDEs with a known structure. Reference \cite{long2017pde} makes a substantial contribution by using neural networks to accurately learn partial differential equations with non constant coefficients.  A neural network is constructed that mimics a forward Euler timestepping scheme and the accuracy of potential models is evaluated based on their future state prediction accuracy.  While seemingly more robust than sparse regression, the method in \cite{long2017pde} does not penalize extraneous terms in the learned PDE and thus falls short of producing optimally parsimonious models.  Furthermore, \cite{long2017pde} only tests the method on a nonlinear problem using a relatively strong Ansatz.  Neural networks were also used in \cite{raissi2017physics1} and \cite{raissi2017physics2} to solve and to estimate parameters in partial differential equations with known terms to a high degree of accuracy.  However, similar to \cite{raissi2017hidden}, it is assumed that the PDE is known up to a set of coefficients.  A more sophisticated neural network approach was used in \cite{raissi2018deep} to learn dynamics of systems with unknown terms.  However, the approach in \cite{raissi2018deep} does not give closed form representations of the dynamics and the resulting neural network model therefore does not give insights into the underlying physics.

In this work, we present a sparse regression framework for identifying PDEs with non-constant coefficients, something that none of the previous PDE discovery methods are equipped to do.  Specifically,  we allow for variation in the value of a coefficient across time or space, but maintain that the active terms in the PDE must be consistent.  This is an important innovation in practice as the parameters of physical systems often vary during the measurement process, so that the parametric dependencies be disambiguated from the PDE itself.  Our method extends the sparse regression frameworks first proposed for PDE discovery~\cite{rudy2017data,schaeffer2017learning} by using group sparsity and results in a more parsimonious and interpretable model than neural networks. We are still limited by the accuracy of numerical differentiation and by the library terms in the sparse regression. Numerical differentiation using neural networks as shown in \cite{raissi2018deep} appears promising as a method for obtaining more accurate time derivatives from noisy data.  The limitation based on terms included in the library seems more permanent.  Any closed-form model expression must be representable with a finite set of building blocks.  Here we only use monomials in our data and its derivatives, since these are the common terms seen in physics, but there is no limitation to the terms included in the library.  

\section{Methods}

The parametric discovery method relies on several foundational mathematical tools. In the following subsections, we will discuss the identification of constant coefficient equations as well as regression methods for group sparsity.  Finally, we combine these ideas to show how one may identify parametric PDEs and suggest a methodology for model selection that balances accuracy and the number of active terms in the PDE.

\subsection{Identification of constant coefficient partial differential equations}

Several recent methods have been proposed for the identification of constant coefficient partial differential equations from data.  In this work we expand on the sparse regression framework, PDE-FIND, used in \cite{rudy2017data}.  We will briefly elaborate on the method and refer the reader to the original paper for details.  The PDE-FIND algorithm provides a principled technique for discovering the underlying PDE from spatial time series measurements alone using a library of candidate functions for the PDE and sparse regression.  For the identification of constant coefficient PDEs, we have a dataset ${\bf U}$, which is a discretization of a function $u(x,t)$ that we assume satisfies the PDE of the form given in (\ref{eq:overview}):
\begin{equation}
u_t= N(u, u_x, u_{xx}, \dots) = \sum_{j=1}^d N_j(u, u_x, u_{xx}, \dots) \xi_j \, .
\label{eq:const_coeff_PDE}
\end{equation}
We assume that the nonlinear expression $N(\cdot )$ may be expanded as a sum of simple monomial basis functions $N_j$ of $u$ and its derivatives.  Note that this sum is not unique and that we can include extra basis functions by simply setting the corresponding $\xi_j$ to be zero.  In PDE-FIND, which constructs an overcomplete library of many possible monomial basis functions and regresses to find $\xi$, sparsity is used to ensure that basis functions that do not appear in the PDE are set to zero in the sum.  

To be more precise, given a dataset ${\bf U} \in \R^{n x m}$ representing $m$ timesteps of a PDE discretized with $n$ gridpoints, we numerically differentiate in both $x$ and $t$ to form the linear regression problem given by (\ref{eq:PDE_FIND})
\begin{equation}
\underbrace{
\begin{pmatrix}
u_t(x_1,t_1) \\
u_t(x_2,t_1) \\
\vdots \\
u_t(x_n,t_m) \\
\end{pmatrix}}_{{\bf U}_t}
=
\underbrace{
\begin{pmatrix}
1 & u(x_1,t_1) & u_x(x_1,t_1) & \hdots & u^3u_{xxx}(x_1,t_1) \\
1 & u(x_2,t_1) & u_x(x_2,t_1) & \hdots & u^3u_{xxx}(x_2,t_1) \\
\vdots & \vdots & \vdots & & \vdots \\
1 & u(x_n,t_m) & u_x(x_n,t_m) & \hdots & u^3u_{xxx}(x_n,t_m) \\
\end{pmatrix}}_{\mathbf{\Theta} ({\bf U})}
\xi\label{eq:PDE_FIND}
\end{equation}
which is a large, overdetermined linear system of equations ${\bf A}{\bf x}={\bf b}$.  Note that here we have shown a problem where derivatives up to third order are multiplied by powers of $u$ up to cubic order, but one could include arbitrarily many library functions.  Solving for $\xi$ and ensuring sparsity gives the PDE.  PDE-FIND has been shown to accurately identify several partial differential equations from data alone.  The sparsity constraint is a regularizer for the linear regression~\cite{rudy2017data}.
\subsection{Group Sparsity}
In a typical sparse regression, we seek a sparse solution to the linear system of equations ${\bf A}{\bf x}={\bf b}$.  Accuracy of the predictor, $\|${\bf A}{\bf x}-{\bf b}$\|$, is balanced against the number of nonzero coefficients in ${\bf x}$.  Thus the sparse regularization enforces a solution ${\bf x}$ with many zeros (which is the variable $\boldsymbol{\xi}$ in (\ref{eq:PDE_FIND})).  In this paper, we use the notion of group sparsity to find time series representing each parameter in the PDE, rather than single values. We group collections of terms in ${\bf x}$ together and seek solutions to ${\bf A}{\bf x}={\bf b}$ that minimize the number of groups with nonzero coefficients.

One well studied method for solving regression problems with group sparsity is the group LASSO (GLASSO)~\cite{friedman2010note} 
\begin{equation}
\label{eq:group_lasso}
\hat{\bf x} =  \underset{\bf w}{\mbox{arg\,min}} \frac{1}{2n} \left\|{\bf b} - \sum_{g \in \mathcal{G}} \mathbf{A}^{(g)} {\bf w}^{(g)}\right\|_2^2 + \lambda \sum_{g \in \mathcal{G}} \|{\bf w}^{(g)}\|_2 .
\end{equation}
Here $\mathcal{G}$ is a collection of groups, each of which contains a subset of the indices enumerating the columns of ${\bf A}$ and coefficients in ${\bf x}$.  Note that the second term in the GLASSO corresponds to a convex relaxation of the number of groups containing a nonzero value.

The concept of group sparsity has been used in several previous methods for identifying dynamics given by ordinary differential equations \cite{chen2017network, schaeffer2017learning_group}.  As it will be shown in the following, we find that the GLASSO performs poorly in the case of identifying PDEs.  We instead use a sequential thresholding method based on ridge regression, similar to the method used in \cite{rudy2017data}, but adapted for group sparsity.  A sequential thresholding method was also used in \cite{schaeffer2017learning_group} for group sparsity but for ordinary and not partial differential equations.  Our method, which we call Sequential Grouped Threshold Ridge Regression (SGTR), is summarized in algorithm \ref{SGTR}.

\begin{center}
\begin{minipage}{0.8\textwidth}
\begin{algorithm}[H]
    \caption{SGTR($\mathbf{A}, {\bf b}, \mathcal{G}, \lambda, \epsilon, \text{maxit}, f({\bf x}) = \|{\bf x} \|_2$)
}
    \label{SGTR}
    \vspace{1 mm}
    
    \# Solves ${\bf x} \approx \mathbf{A}^{-1}{\bf b}$ with sparsity imposed on groups in $\mathcal{G}$\\
    
    \vspace{1 mm}
    \# Initialize coefficients with ridge regression\\
    ${\bf x} = \mbox{arg\,min}_{\bf w} \|{\bf b} - \mathbf{A}{\bf w}\|_2^2 + \lambda \|{\bf w} \|_2^2$
    
    \vspace{3 mm}
    \# Threshhold groups with small $f$ and repeat
    
    for $iter = 1,\hdots,\text{maxit}$:\\
    
	\hspace{1 cm} \# Remove groups with sufficiently small $f({\bf x}^{(g)})$
	
	\hspace{1 cm} $\mathcal{G} = \{g \in \mathcal{G} : f({\bf x}^{(g)}) > \epsilon\}$\\

	\hspace{1 cm} \# Refit these groups (note this sets ${\bf x}^{(g)}=0$ for $g \not\in \mathcal{G}$)

	\hspace{1 cm} ${\bf x} = \mbox{arg\,min}_{\bf w} \|{\bf b} - \sum_{g \in \mathcal{G}}\mathbf{A}^{(g)} {\bf w}^{(g)}\|_2^2 + \lambda \|{\bf w}\|_2^2$
	\vspace{3 mm}
	
	\# Get unbiased estimates of coefficients after finding sparsity\\
	${\bf x}^{(\mathcal{G})} = \mbox{arg\,min}_{\bf w} \|{\bf b} - \sum_{g \in \mathcal{G}}\mathbf{A}^{(g)} {\bf w}^{(g)}\|_2^2$\\
	
	return ${\bf x}$
	
\end{algorithm}
\end{minipage}

\end{center}

\vspace{3 mm}

Throughout the training, $\mathcal{G}$ tracks the groups that have nonzero coefficients, and it is paired down as we threshold coefficients with sufficiently small relevance, as measured by $f$.  We use the 2-norm of the coefficients in each group for $f$ but one could also consider arbitrary functions.  In particular, for problems where the coefficients within each group have a natural ordering, as they do in our case as time series or spatial functions, one could consider smoothness or other properties of the functions.  In practice, we normalize each column of ${\bf A}$ and ${\bf b}$ so that differences in scale between the groups do not affect the result of the algorithm.  For the GLASSO we always perform an unregularized least squares regression on the nonzero coefficients after the sparsity pattern has been discovered to debias the coefficients.  We found SGTR to outperform the GLASSO for the problem of correctly identifying the active terms in parametric PDEs.

\subsection{Data-driven identification of parametric partial differential equations}

In the identification of parametric PDEs, we consider equations of the form
\begin{equation}
\label{eq:parametric_PDE}
u_t = N(u, u_x, \hdots, \mu(t)) = \sum_{j=1}^D N_j(u, u_x, \hdots) \xi_j (t).
\end{equation}
Note that this equation is similar to (\ref{eq:const_coeff_PDE}) but has time-varying parametric dependence.  To capture spatial variation in the coefficients, we simply replace $\xi(t)$ with $\xi(x)$.  The PDE is assumed to contain a small number of active terms, each with a time varying coefficient $\xi (t)$.  We seek to solve two problems:  (i) determine which coefficients are nonzero and (ii) find the values of the coefficients for each $\xi_j$ at each timestep or spatial location for which we have data.

For time dependent problems, we construct a separate regression for each timestep, allowing for variation in the PDE between timesteps.  Similar to the PDE-FIND method, we construct a library of candidate functions for the PDE using monomials in our data and its derivatives so that
\begin{equation}
\label{eq:single_timestep_library}
\mathbf{\Theta} \left(u^{(j)}\right) = \begin{pmatrix}
\vline & \vline & \vline & \vline \\
1 & u^{(j)} & \hdots & u^3u^{(j)}_{xxx} \\
\vline & \vline & \vline & \vline 
\end{pmatrix}
\end{equation}
where the set of $m$ equations is given by  
\begin{equation}
\label{eq:parametric_PDE_FIND_setup_2}
u_t^{(j)} = \mathbf{\Theta} \left(u^{(j)}\right)\xi^{(j)} ,\,\, j = 1,\hdots, m \, .
\end{equation}
Our goal is to solve the set of equations given by (\ref{eq:parametric_PDE_FIND_setup_2}) with the constraint that each $\xi^{(j)}$ is sparse and that they all share the same sparsity pattern.  That is, we want a fixed set of active terms in the PDE.  To do this, we consider the set of equations as a single linear system and use group sparsity.  Expressing the system of equations for the parametric equation as a single linear system we get the block diagonal structure given by 
\begin{equation}
\label{eq:parametric_PDE_FIND_setup}
\begin{pmatrix}
u_t^{(1)}\\
u_t^{(2)}\\
\vdots \\
u_t^{(m)}
\end{pmatrix} = 
\underbrace{\begin{pmatrix}
\mathbf{\Theta} \left(u^{(1)}\right) \\
& \mathbf{\Theta} \left(u^{(2)}\right) \\
&& \ddots \\
&&& \mathbf{\Theta} \left(u^{(m)}\right)
\end{pmatrix}}_{\tilde{\mathbf{\Theta}}}
\begin{pmatrix}
\xi^{(1)}\\
\xi^{(2)}\\
\vdots \\
\xi^{(m)} \, 
\end{pmatrix}.
\end{equation}

We solve (\ref{eq:parametric_PDE_FIND_setup}) using SGTR with columns of the block diagonal library matrix grouped by their corresponding term in the PDE.  Thus for $m$ timesteps and $d$ candidate functions in the library, groups are defined as $\mathcal{G} = \{ {j + d\cdot i : i = 1,\dots, m} : j = 1,\hdots, d\}$.  This ensures a sparse solution to the PDE while also allowing arbitrary time series for each variable.  To obtain the correct level of sparsity, we separate 20\% of the data from each timestep to use as a validation set and search over the parameter $\lambda$ in the SGTR algorithm using cross validation to find the optimal $\lambda$.  For problems with spatial, rather than temporal, variation in the coefficients, we simply group by spatial rather than time coordinate.  A similar block diagonal structure is obtained but with $n$ blocks of size $m \times d$ rather than $m$ blocks of size $n \times d$.  The groups are defined by $\mathcal{G} = \{ {j + d\cdot i : i = 1,\dots, n} : j = 1,\hdots, d\}$.

Since we are evaluating the relevance of groups based on their norm, it is important to consider differences in the scale of the candidate functions.  For example, if $u \sim \mathcal{O}(10^{-2})$ then a cubic function will be $\mathcal{O}(10^{-6})$ and relatively large coefficients multiplying this data may not have a large effect on the dynamics, but will not be removed by the hard threshold due to its size.  To remedy this, we normalize each candidate function represented in $\mathbf{\Theta}$ as well as each $u_t^{(j)}$ to have unit length prior to the group thresholding algorithm and then correct for the normalization after we have discovered the correct sparsity pattern.

\subsection{Model Selection}

For each model we test both the GLASSO as well as SGTR using an exhaustive range of parameter values.  Let $\tilde{\mathbf{\Theta}}$ denote the block diagonal matrix $\mathbf{\Theta}$ shown in equation \eqref{eq:parametric_PDE_FIND_setup} but with all columns normalized to have unit length and $\tilde{u}_t$ be the vector of all time derivatives having been normalized to unit length so that $\|\tilde{u}_t\| = \sqrt{m}$.  For the GLASSO, we find the minimal value of $\lambda$ that will set all coefficients to zero which is given by 
\begin{equation}
\label{eq:lam_max}
\lambda_{\mbox{max}} = \underset{g \in \mathcal{G}}{\mbox{max}} \,\frac{1}{n}\|\tilde{\mathbf{\Theta}}^{(g)^T}\tilde{u}_t\|_2 .
\end{equation}
We check 50 evenly spaced  values of $\lambda$ between $10^{-5} \lambda_{\mbox{max}}$ and $\lambda_{\mbox{max}}$ on a logarithmic scale.

For SGTR, we search over the range of tolerances between $\epsilon_{min}$ and $\epsilon_{max}$ defined as
\begin{equation}
\label{eq:max_min_tol}
\epsilon_{max/min} = \underset{g \in \mathcal{G}}{max/min}\, \|\xi_{\text{ridge}}^{(g)}\|_2,
\end{equation}
where $\xi_{\text{ridge}} = (\tilde{\mathbf{\Theta}}^T\tilde{\mathbf{\Theta}} + \lambda I)^{-1}\tilde{\mathbf{\Theta}}^T\tilde{u}_t$.  Note that by definition, $\epsilon_{min}$ is the minimum tolerance that has any effect on the sparsity of the predictor and $\epsilon_{max}$ is the minimum tolerance that guarantees all coefficients to be zero.  A set of 50 intermediate tolerances, equally spaced on a logarithmic scale, is tested between $\epsilon_{min}$ and $\epsilon_{max}$.

To select the optimal model generated via each method, we evaluate the models using the AIC-inspired loss function 
\begin{equation}
\label{eq:AIC}
\mathcal{L}(\xi) = N \ln\left(\dfrac{\|\tilde{\mathbf{\Theta}}\xi-\tilde{u}_t\|_2^2}{N} + \epsilon \right)  + 2k
\end{equation}
where $k$ is the number of nonzero coefficients in the identified PDE, $\|\xi\|_0/m$, and $N$ is the number of rows in $\mathbf{\Theta}$, which is equal to the size of our original dataset $u$.


Equation \eqref{eq:AIC} is closely related to the Akaike Information Criterion (AIC) \cite{akaike1974new}.  Typically, the mean square error of a linear model is used to evaluate goodness of fit, but in our case there is error in computing the time derivative $u_t$, so we assume that any linear model which perfectly fits the data is overfit.  We have added $\epsilon = 10^{-5}$ to the mean square error of each model as a floor in order to avoid overfitting.  Without this addition, our algorithm selects insufficiently parsimonious representations of the dynamics.

\begin{figure}[t]
\centering
\includegraphics[width=1.0\textwidth]{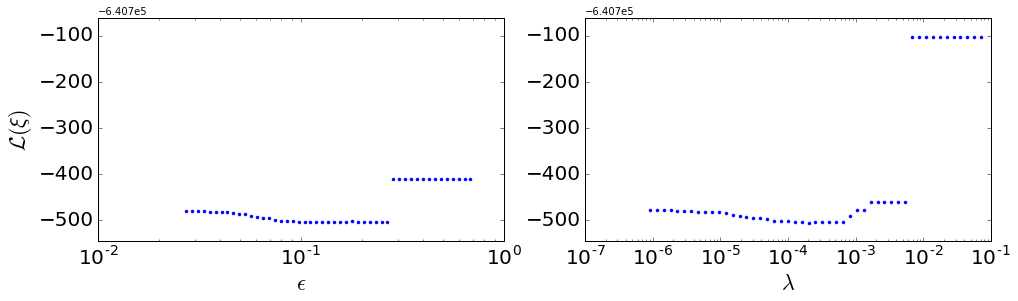}
\caption{Example of loss function evaluated for a number of candidate models for the parametric Burgers' equation.  Library included derivatives up to 4th order muliplying powers of $u$ up to cubic.  Left: 50 models obtained via SGTR algorithm using values of $\epsilon$ between $\epsilon_{min}$ and $\epsilon_{max}$.  Right: 50 models obtained via GLASSO for $\lambda$ between $10^{-5}\lambda_{max}$ and $\lambda_{max}$.}
\label{fig:aic_fig}
\end{figure}

Figure~\ref{fig:aic_fig} illustrates the loss function from equation \eqref{eq:AIC} evaluated on models derived from 50 values of $\epsilon$ and $\lambda$ using SGTR and GLASSO respectively.  Initially, a low penalty in each algorithm yields a model that is overfit to the data given our sparsity criteria.  For an intermediate value of $\epsilon$ or $\lambda$, a more parsimonious but still predictive model is obtained.  For sufficiently high values, the model is too sparse and is no longer predictive.

\section{Computational Results of Parametric PDE Discovery}

%
%

We test our method for the discovery of parametric partial differential equations on four canonical models; Burgers' equation with a time varying nonlinear term, the Navier-Stokes equation for vorticity with a jump in Reynolds number, a spatially dependent advection equation, and a spatially dependent Kuramoto-Sivashinsky equation. In each case the method is also tested after introducing white noise with mean magnitude equal to 1\% of the $L^2$-norm \ale{$L^2$ or $\ell_2$}of the dataset. The method is able to accurately identify the dynamics in each case except the Kuramoto-Sivashinsky equation, where the inclusion of a fourth order derivative makes numerical evaluation with noise highly challenging.  A comparison with GLASSO regression is also given for a number of the examples.

\subsection{Burgers' Equation with Diffusive Regularization}

\begin{figure}[t]
\centering
\includegraphics[width=1.0\textwidth]{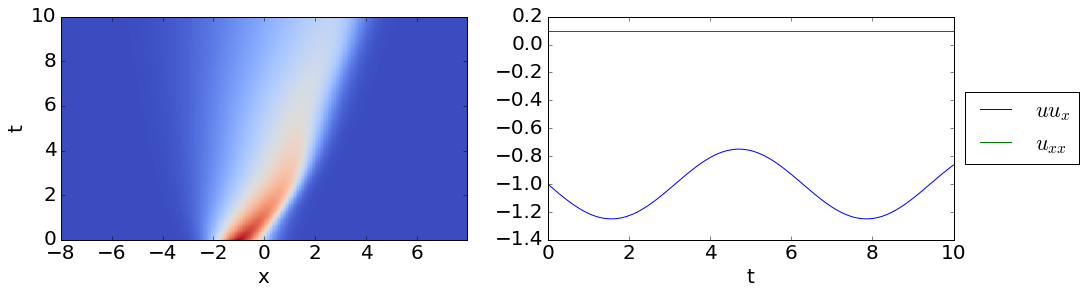}
\caption{Left: dataset for identification of the parametric diffusive Burgers' equation.  Here the PDE was evolved on the interval $[-8,8]$ with periodic boundary conditions and $t \in [0,10]$.  Right:  Coefficients for the terms in the parametric Burgers' equation.  The diffusion was held constant at 0.1 while the nonlinear advection as coefficient is given by $a(t)=-(1+sin(t)/4)$.}
\label{fig:parametric_burgers_solution}
\end{figure}

To test the parametric discovery of PDEs, we consider a solution of Burgers' equation with a sinusoidally oscillating coefficient $a(t)$ for the nonlinear advection term
\begin{equation}
\label{eq:burgers}
\begin{aligned} u_t &= a(t) uu_x + 0.1 u_{xx}\\[.1in]
a(t) &= -\left(1 + \frac{\sin(t)}{4}\right)
\end{aligned}
\end{equation}
where a small amount of diffusion is added to regularize the evolution dynamics.

The time dependent Burgers' equation was solved numerically using a spectral method on the interval $[-8,8]$ with periodic boundary conditions and $t \in [0,10]$ with $n = 256$ grid points and $m = 256$ time steps.  We search for parsimonious representations of the dynamics by including powers of $u$ up to cubic order, which can be multiplied by derivatives of $u$ up to fourth order.  For the noise-free dataset we use the discrete Fourier transform for computing derivatives.  For the noisy dataset, we use polynomial interpolation to smooth the derivatives~\cite{rudy2017data}.

\begin{figure}[t]
\centering
\includegraphics[width=1.0\textwidth]{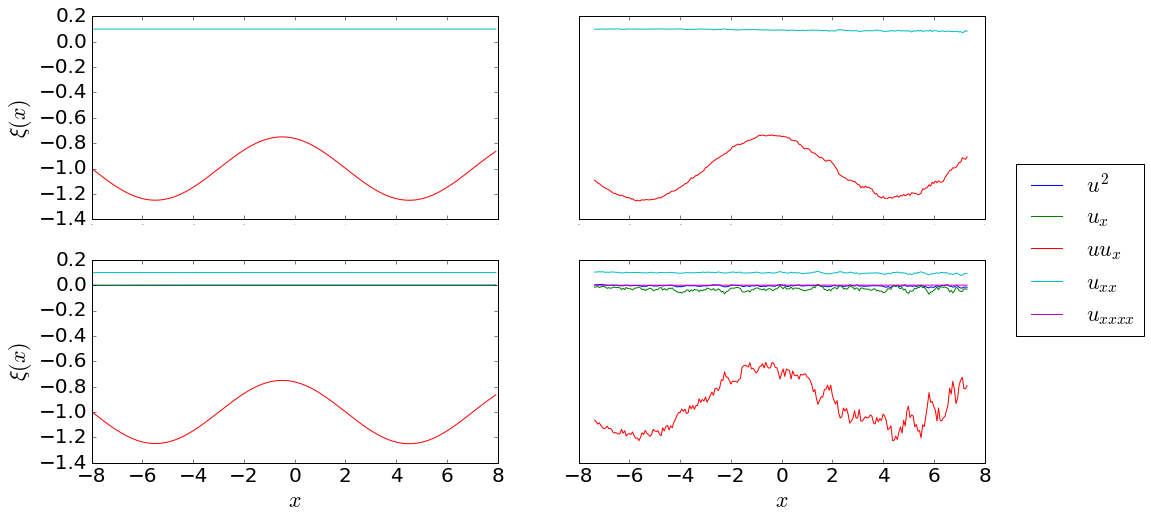}
\caption{Time series discovered for the coefficients of the parametric Burgers' equation.  Top row: SGTR method, which correctly identifies the two terms.  Bottom row:  GLASSO method which adds several additional (incorrect) terms to the model.   The left panels are noise-free, while the right panels contain 1\% noise.  This parametric dependency is illustrated in Fig.~\ref{fig:param}b.}
\label{fig:parametric_burgers}
\end{figure}

The resulting time series for the identified nonzero coefficients are shown in Fig.~\ref{fig:parametric_burgers}.  SGTR correctly identified the active terms in the PDE for both the noise-free and noisy datasets, whereas GLASSO fails in both cases to produce the correct PDE and its parametric dependencies.

\subsection{Navier-Stokes:  Flow around a cylinder}
\begin{figure}[t]
\centering
\includegraphics[width=1.0\textwidth]{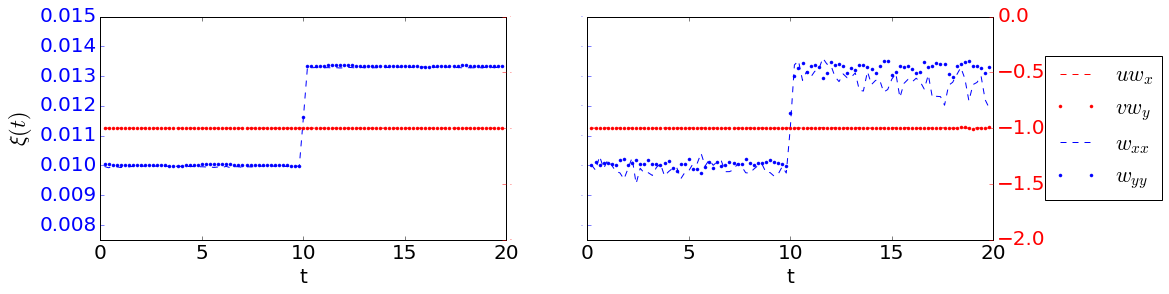}
\caption{Identified time series for coefficients for the Navier Stokes equation.  Distinct axes are used to highlight jump in Reynolds number.  Left: no noise.  Right: 1\% noise}
\label{fig:parametric_ns}
\end{figure}

We consider the fluid flow around a circular cylinder by simulating the Navier-Stokes vorticity equation
\begin{equation}
\label{eq:navier_stokes}
\omega_t + \mathbf{u}\cdot \nabla \omega = \dfrac{1}{\nu (t)} \Delta \omega.
\end{equation}
Data is generated using the Immersed Boundary Projection Method (IBPM) \cite{taira:07ibfs,taira:fastIBPM} with $n_x = 449$ and $n_y=199$ spatial points in $x$ and $y$ respectively, and 1000 timesteps with $dt = 0.02$.  The Reynolds number is adjusted half way through the simulation from $\nu = 100$ initially to $\nu = 75$.  This is representative of the fluid velocity exhibiting a sudden decrease midway through the data collection.  Our library of candidate functions is constructed using up to second order derivatives of the vorticity and multiplying by up to quadratic functions of the data.  To keep the size of the machine learning problem tractable, we subsample 1000 random spatial location from the wake of the cylinder to construct our library at every tenth timestep~\cite{rudy2017data}.  For the noise-free dataset, far fewer points are needed to accurately identify the dynamics.  We suspect that with a more careful treatment of the numerical differentiation in the case of noisy data, such as that used in \cite{raissi2018deep}, the same would be true for the dataset with artificial noise, however such work is not the focus of this paper.
The identified time series for the Navier-Stokes equation are shown in Fig.~\ref{fig:parametric_ns}.  SGTR and GLASSO both correctly identify the active terms in the PDE.

\begin{figure}[t]
\centering
\includegraphics[width=1.0\textwidth]{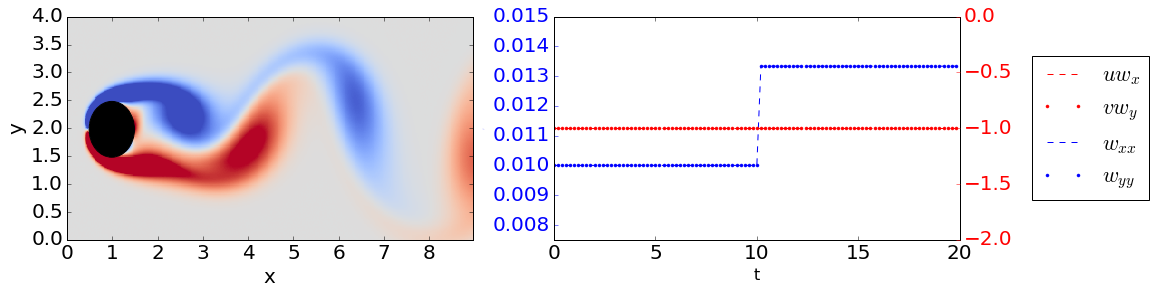}
\caption{Left: dataset for identification of the parametric Navier-Stokes equation (\ref{eq:navier_stokes}).   Right: coefficients for Navier-Stokes equations exhibiting jump in Reynolds number from 100 to 75 at $t=10$.  This parametric dependency is illustrated in Fig.~\ref{fig:param}a.}
\label{fig:parametric_ns_solution}
\end{figure}

\subsection{Spatially Dependent Advection-Diffusion Equation}
%

\begin{figure}[t]
\centering
\includegraphics[width=1.0\textwidth]{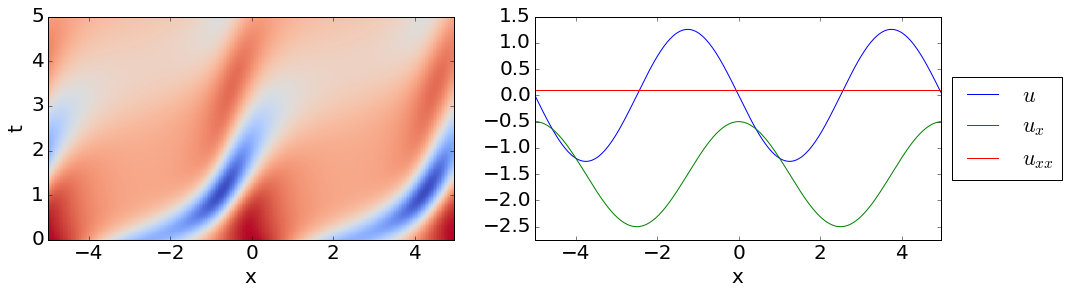}
\caption{Left: dataset for identification of the spatially dependent advection diffusion equation.  Right:  Spatial dependence of PDE.  In this case, the loadings $\xi_j(t)$ in (\ref{eq:parametric_PDE}) are replaced by $\xi_j(x)$.}
\label{fig:advection_solution}
\end{figure}
\begin{figure}[h]
\centering
\includegraphics[width=1.0\textwidth]{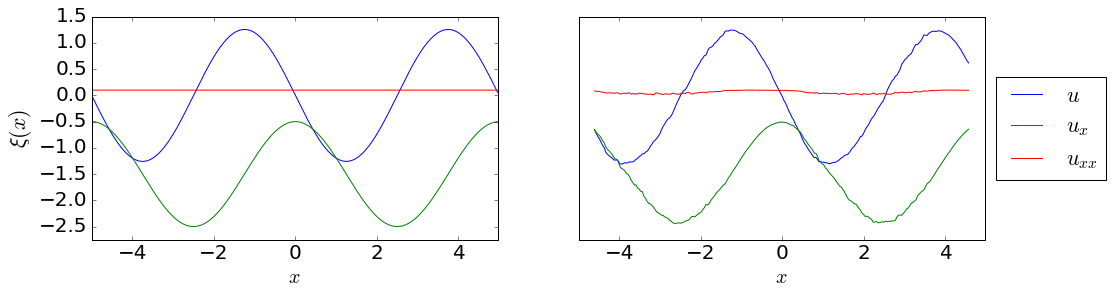}
\caption{Spatial dependence of advection diffusion equation.  Left: no noise.  Right: 1\% noise.  Both SGTR and GLASSO correctly identified the active terms.}
\label{fig:spatial_advection}
\end{figure}

The advection-diffusion equation is a simple model for the transport of a physical quantity in a velocity field with diffusion.  Here, we adapt the equation to have a spatially dependent velocity
\begin{equation}
\label{eq:advection}
u_t = (c(x) u)_x + \epsilon u_{xx} = c(x) u_x + c'(x)u+ \epsilon u_{xx}
\end{equation}
which models transport through a spatially varying vector field due to $c=c(x)$.
The PDE is solved on a periodic domain $[-L,L]$ with $L=5$, $\epsilon = 0.1$, and $c(x) = -1.5 + \cos(2\pi x /L)$ using a spectral method with $n = 256$ and $m = 256$.  The library consists of powers of $u$ up to cubic, multiplied by derivatives of $u$ up to fourth order.
Results for the advection-diffusion equation are shown in Fig.~\ref{fig:spatial_advection}.  In the noise-free and noisy datasets, both SGTR and GLASSO Correctly identify the active terms in the PDE.

\subsection{Spatially Dependent Kuramoto-Sivashinsky Equation}
\begin{figure}[t]
\centering
\includegraphics[width=1.0\textwidth]{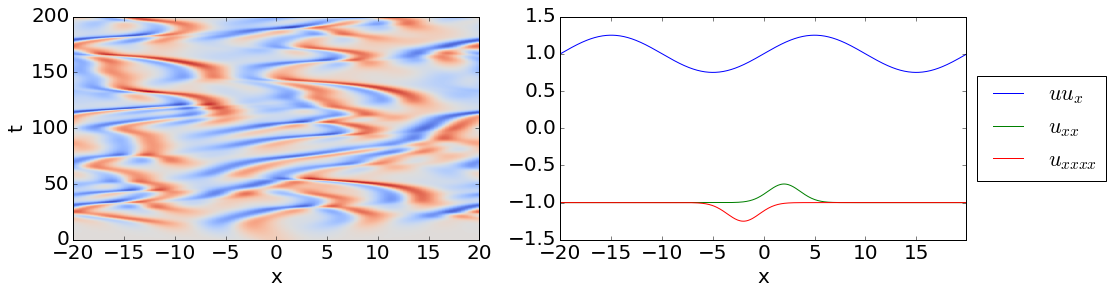}
\vspace{-.2in}
\caption{Left:  dataset for identification of the spatially dependent Kuramoto Sivashinsky equation. Right: parametric dependency of the governing equations.}
\label{fig:ks_solution}
\end{figure}
We now test the method on a Kuramoto-Sivashinsky equation with spatially varying coefficients 
\begin{equation}
\label{eq:KS}
u_t = a(x) uu_x + b(x) u_{xx} + c(x) u_{xxxx} .
\end{equation}
We use a periodic domain $[-L,L]$ with $L=20$ and coefficients $a(x) = 1 + \sin(2\pi x /L) /4$, $b(x) = -1 + e^{-(x-2)^2/5}/4$ and $c(x) =  -1 - e^{-(x+2)^2/5}/4$.

\begin{figure}[t]
\centering
\includegraphics[width=1.0\textwidth]{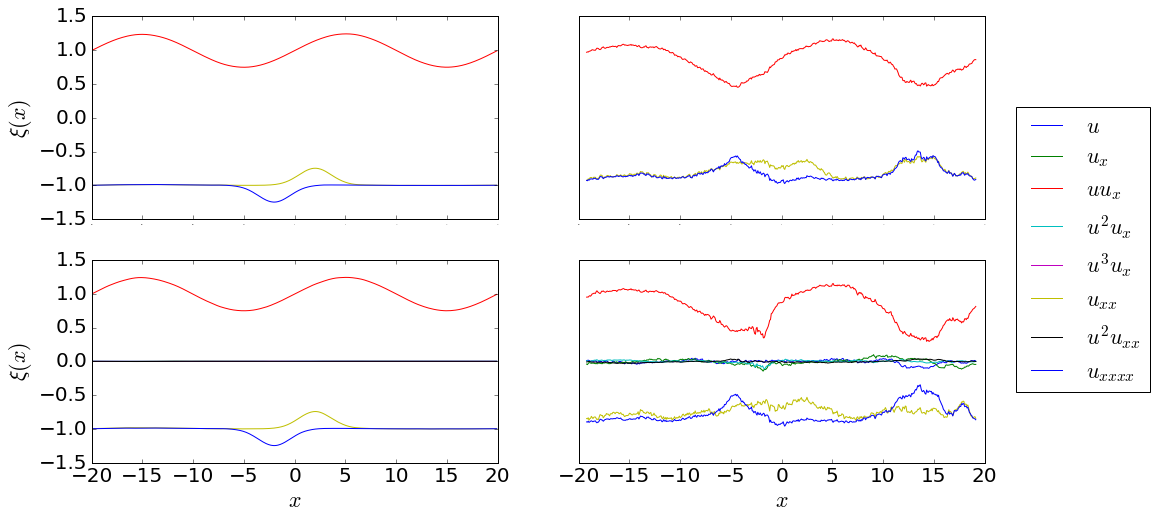}
\vspace{-.2in}
\caption{Spatial dependence of Kuramoto-Sivashinsky.  top row: SGTR.  Bottom row: GLASSO.  Left: no noise. Right 0.01\% noise using SGTR.  SGTR detects correct sparsity with significant parameter error.  GLASSO does not correctly identify the parsimonious model, nor does it do a good job at predicting the correct parametric values.}
\label{fig:spatial_ks}
\end{figure}

The equation is solved numerically to $t = 200$ using $n=512$ grid points and $m=1024$ timesteps.  The second half of the data set is used, so as to only consider the region where the dynamics exhibited spatio-temporal chaos, resulting in a dataset containing 512 snapshots of 512 gridpoints.  Since the Kuramoto Sivashinsky equation involves a fourth order derivative, it is very difficult to correctly identify it with noisy data since it is exceptionally difficult to accurately compute the fourth derivative.  Indeed, our method fails to correctly identify the active terms when 1\% noise is added.  With 0.01\% noise the correct terms were identified but with substantial error in coefficient value.  We suspect that this shortcoming could at least be partially remedied by a more careful treatment of the numerical differentiation such as in \cite{bruno2012numerical}.  The results of our parametric identification are shown in Fig.~\ref{fig:spatial_ks}.

\section{Discussion}

We have presented a method for identifying governing laws for physical systems which exhibit either spatially or temporally dependent behavior.  The method builds on a growing body of work in the applied mathematics and machine learning community that seeks to automate the process of discovering physical laws.  To the best of our knowledge, our method is the first approach for deriving parsimonious PDE expressions of spatio-temporal system in the case of non-constant coefficients.  Specifically, we can disambiguate between the governing PDE evolution and its parametric dependencies.  In all examples, the SGTR algorithm outperformed the GLASSO in correctly identifying active terms in the PDE.  Errors from the latter were generally in the form of extra terms added with small coefficient values throughout the time series.  It may seem reasonable to threshold these time series after the discovery algorithm, but doing so assumes that the terms importance in the PDE is directly related to its magnitude, an assumption which we do not make given the normalization prior to sparse regression.

In this work we have split the data into distinct timesteps or spatial locations in order to find PDE models for each subset of the data, resulting in coefficients that can vary in space or time.  However, with a sufficiently fine grid, it seems feasible that one could bin the data by areas localized in space and time to determine a coefficient varying in both space and time with some loss of resolution.  This same result may be achievable in a more stable manner by introducing a sparsity term to the work in \cite{long2017pde}.

As is the case with other sparse regression methods for identifying dynamical systems, this method is constrained by the ability of the user to accurately differentiate data.  For ordinary differential equations, this may be circumvented by looking at the weak form of the dynamics \cite{schaeffer2017extracting}, but doing so for PDEs seems difficult since there are derivatives that need to be evaluated with respect to multiple variables.  We find the automatic differentiation approach used in \cite{raissi2018deep} promising and suspect that the inclusion of neural network based differentiation could radically improve the ability of our method to identify dynamics from noisy data.  With sufficient knowledge of data it may also be possible to obtain better estimates through tuning the polynomial based differentiation \cite{bruno2012numerical}.

Automating the identification of closed form physical laws from data will hopefully boost scientific progress in areas where deriving the same laws from first principals proves intractable.  There are several limitations to many methods proposed in the field thus far.  In particular, current methods have generally studied equations of the form $u_t = N(u,x,t)$ but many equations in physics are not in this class.   Indeed, if measuring a system with parametric dependencies, then past methods are be unable to disambiguate between the evolution dynamics and its parametric dependencies $\mu(t)$, thus greatly limiting model discovery.  
There is also a trade off between methods that are able to derive parsimonious representations, though which are limited to a finite set of library elements, and those that use black box models to represent larger classes of possible functions.  The researcher may also find difficulties in attempting to infer dynamics from the wrong set of measurements.  For example, one could not derive the Schr\"{o}dinger by only looking at measurements of intensity.  While not addressing these issues, this work makes a step towards generalizing the class of equations which may be accurately identified via machine learning methods.\\

\noindent\footnotesize{Code: \href{https://github.com/snagcliffs/parametric-discovery}{https://github.com/snagcliffs/parametric-discovery}
}

\begin{appendix} 

\end{appendix}


\bibliographystyle{siamplain}
\bibliography{ROM_FIND_bibliography}

\end{document}